\documentclass[a4paper,12pt]{amsart}
\usepackage{amsmath,amssymb,amsthm}
\usepackage{multirow}
\usepackage{longtable}

\usepackage{url}
\usepackage{comment}

\newtheorem*{Thm}{Theorem}

\newcommand{\Q}[1]{\mathbb{Q}(\sqrt{#1})}
\newcommand{\Z}{\mathbb{Z}}

\newcommand{\calO}{\mathcal{O}}

\newcommand{\qf}[1]{\langle #1 \rangle}
\newcommand{\conj}[1]{\overline{#1}}
\newcommand{\binlattice}[4]{\begin{bmatrix}
  #1 & #2 \\
  #3 & #4
\end{bmatrix}}

\newcommand{\comega}{{\conj\omega}}
\newcommand{\nequiv}{\not\equiv}
\newcommand{\Tr}{\operatorname{Tr}}

\newcommand{\rank}{\operatorname{rank}}

\begin{document}

\title{Simple proofs for universal binary Hermitian lattices}

\author[Poo-Sung Park]{Poo-Sung Park}
\address{School of Computational Sciences, Korea Institute institute for
Advanced Study, Hoegiro 87, Dongdaemun-gu, Seoul, 130-722, Korea}
\email{sung@kias.re.kr}



%\begin{abstract}
%If a positive definite Hermitian lattice represents all positive
%integers, we call it universal. Several mathematicians, including
%the author, found 25 universal binary Hermitian lattices. But their
%ad hoc proofs are complicated. We give simple and unified proofs.
%\end{abstract}

\maketitle

\section{Introduction}
%%%%%%%%%%%%%%%%%%%%%%%%%%%%%%%%%%%%%%%%%%%%%%%%%%%%%%%%%%%%%%%%%%%%%%%%%%%

It has been a central problem in the theory of quadratic forms to
find integers represented by quadratic forms. The celebrated
\emph{Four Square Theorem} by Lagrange \cite{Lagrange} was an
outstanding result in this study. Ramanujan generalized this theorem
and found $54$ positive definite quaternary quadratic forms which
represent all positive integers \cite{Ramanujan}. We call a positive
definite quadratic form \emph{universal}, if it represents all
positive integers. The classification of nondiagonal universal
classical quadratic forms was completed by Conway and Schneeberger
using their \emph{Fifteen Theorem} in 2000 \cite{Conway},
\cite{Bhargava}. The theorem states that if a positive definite
classical quadratic form (with four or more variables) represents up
to $15$, it is universal.

In 1997 Earnest and Khosravani defined universal Hermitian forms and
they sought 13 positive definite \emph{binary} Hermitian forms over
imaginary quadratic fields of class number one
\cite{Earnest-Khosravani}. Iwabuchi extended the result to imaginary
quadratic fields of class number bigger than one and he found 9
binary Hermitian lattices (as a generalization of Hermitian forms)
\cite{Iwabuchi}. Jae-Heon Kim and the author complete the list by
appending 3 universal binary Hermitian forms \cite{Kim-Park}.
Moreover, Kim, Kim and the author found an analogous result to
Fifteen Theorem: If a positive definite Hermitian lattice represents
up to $15$, then it represents all positive integers
\cite{Kim-Kim-Park}. The proof was more complicated than that of the
Conway-Scheeberger Theorem for it contains nonclassical quadratic
forms. The criterion, 290-Theorem, for universal nonclassical
quadratic forms was recently proved by Bhargava and Hanke
\cite{Bhargava-Hanke}.

In the present article we give simple and unified proofs for
universal binary Hermitian lattices. Although the three papers
(\cite{Earnest-Khosravani}, \cite{Iwabuchi}, \cite{Kim-Park})
proposed proofs, they were complicated and used local properties of
Hermitian forms. But, here, we use merely well known results about
quadratic forms.

\section{Notations}

Let $E = \Q{-m}$ for a positive squarefree integer $m$. Denote the
$\mathbb{Q}$-involution by $\conj{\,\cdot\,}$ and the ring of
integers by $\calO = \calO_E$. The generators of $\calO$ are $1$ and
$\omega$ over $\Z$ where $\omega = \sqrt{-m}$ if $m \nequiv 3
\pmod{4}$ and $\omega = \frac{1+\sqrt{-m}}2$ otherwise.

A finitely generated $\calO$-module $L$ is called a \emph{Hermitian
lattice} over $\calO$ if there exists a nondegerate Hermitian space
$(E \otimes_{\calO} L, H)$ over $E$. We consider only positive
definite integral lattices. That is, we assume that
$H(\mathbf{x},\mathbf{y}) \in \calO$ for all $\mathbf{x},\mathbf{y}
\in L$ and $H(\mathbf{x}):=H(\mathbf{x},\mathbf{x}) > 0$ if $0 \ne
\mathbf{x} \in L$.

If the class number of $E$ is one, all Hermitian lattices $L$ are
free and we can write
\[
    L = \calO\mathbf{v}_1 + \calO\mathbf{v}_2 + \dotsb + \calO\mathbf{v}_n
\]
where $n = \rank L = \dim_E{E \otimes L}$. Then the Gram matrix for
$L$ is defined as $M_L = [H(\mathbf{v}_i,\mathbf{v}_j)]_{1 \leq i,j
\leq n}$.

A nonfree Hermitian lattice $L$ can be written as
\[
L = \calO\mathbf{v}_1 + \calO\mathbf{v}_2 + \dotsb +
\calO\mathbf{v}_{n-1} + \mathcal{A}\mathbf{v}_n
\]
with a nonprincipal ideal $\mathcal{A}$ in $\calO$ \cite{OM}. Since
$\mathcal{A}$ is generated by two elements $\alpha, \beta \in
\calO$, we can rewrite
\[
L = \calO\mathbf{v}_1 + \calO\mathbf{v}_2 + \dotsb +
\calO\mathbf{v}_{n-1} + \calO(\alpha\mathbf{v}_n)  +
\calO(\beta\mathbf{v}_n).
\]
So we can deal with $L$ as if $L$ were a \emph{free} Hermitian
lattice of rank $n+1$. The (formal) Gram matrix is defined as an
$(n+1) \times (n+1)$-matrix
\[
    M_L = \begin{bmatrix}
        H(\mathbf{v}_1, \mathbf{v}_1)
            & \cdots
                & H(\mathbf{v}_1, \mathbf{v}_{n-1})
                    & H(\mathbf{v}_1, \alpha \mathbf{v}_n)
                        & H(\mathbf{v}_1, \beta \mathbf{v}_n) \\
        \vdots
            & \ddots
                & \vdots
                    & \vdots
                        & \vdots \\
        H(\alpha \mathbf{v}_n, \mathbf{v}_1)
            & \cdots
                & H(\alpha \mathbf{v}_n, \mathbf{v}_{n-1})
                    & H(\alpha \mathbf{v}_n, \alpha \mathbf{v}_n)
                        & H(\alpha \mathbf{v}_n, \beta \mathbf{v}_n) \\
        H(\beta \mathbf{v}_n, \mathbf{v}_1)
            & \cdots
                & H(\beta \mathbf{v}_n, \mathbf{v}_{n-1})
                    & H(\beta \mathbf{v}_n, \alpha \mathbf{v}_n)
                        & H(\beta \mathbf{v}_n, \beta \mathbf{v}_n) \\
    \end{bmatrix}
\]
whose rank is still $n$, though.

If $L = L_1 \oplus L_2$ and $H(L_1,L_2) = \{ 0 \}$, then we write $L
= L_1 \perp L_2$. If $L$ is a Hermitian lattice generated by only
one vector $\mathbf{v}$, then we write $L = \qf{H(\mathbf{v})}$.
Also $\qf{H(\mathbf{v}_1)} \perp \dotsb \perp \qf{H(\mathbf{v}_n)}$
is written as $\qf{H(\mathbf{v}_1), \dotsc, H(\mathbf{v}_n)}$. From
now on we identify a Hermitian lattice $L$ and its (formal) Gram
matrix $M_L$.

\section{Main Result}

Earnest, Khosravani, Iwabuchi, Kim and the author found all
universal binary Hermitian lattices over imaginary quadratic fields.

\begin{Thm}
There are 25 universal binary Hermitian lattices over imaginary
quadratic fields up to isometry. \small
\begin{align*}
\Q{-1}:\, &\qf{1,1}, \qf{1,2}, \qf{1,3} \\
\Q{-2}:\, &\qf{1,1}, \qf{1,2}, \qf{1,3}, \qf{1,4}, \qf{1,5} \\
\Q{-3}:\, &\qf{1,1}, \qf{1,2} \\
\Q{-5}:\, &\qf{1,2}, \qf{1}\perp\binlattice{2}{-1+\omega_{5}}{-1+\comega_{5}}3 \\
\Q{-6}:\, &\qf{1}\perp\binlattice{2}{\omega_6}{\comega_6}3 \\
\Q{-7}:\, &\qf{1,1}, \qf{1,2}, \qf{1,3} \\
\Q{-10}:\, &\qf{1}\perp\binlattice{2}{\omega_{10}}{\comega_{10}}5 \\
\Q{-11}:\, &\qf{1,1}, \qf{1,2} \\
\Q{-15}:\, &\qf{1}\perp\binlattice{2}{\omega_{15}}{\comega_{15}}2 \\
\Q{-19}:\, &\qf{1,2} \\
\Q{-23}:\, &\qf{1}\perp\binlattice{2}{\omega_{23}}{\comega_{23}}3,
\qf{1}\perp\binlattice{2}{-1+\omega_{23}}{-1+\comega_{23}}3 \\
\Q{-31}:\, &\qf{1}\perp\binlattice{2}{\omega_{31}}{\comega_{31}}4,
\qf{1}\perp\binlattice{2}{-1+\omega_{31}}{-1+\comega_{31}}4
\end{align*}
\end{Thm}

We can associate an $n$-dimensional Hermitian space $(V,H)$ over $E$
with a $2n$-dimensional quadratic space $(\widetilde{V},B_H)$ over
$\mathbb{Q}$ by considering $V$ as a vector space over $\mathbb{Q}$
and defining a bilinear map $B_H(\mathbf{x},\mathbf{y}) =
\frac12\Tr_{E/\mathbb{Q}}H(\mathbf{x},\mathbf{y})$. Thus, to prove
the universality of a given Hermitian lattice, we may show that the
associated quadratic form represents all positive integers.

For $m \nequiv 3 \pmod{4}$ the quadratic forms associated to
$\emph{free}$ Hermitian lattices are diagonal. So their
universalities can be checked by Ramanujan's list. The quadratic
forms associated to
$\qf{1}\perp\binlattice{2}{\omega_6}{\comega_6}{3}$ and
$\qf{1}\perp\binlattice{2}{\omega_{10}}{\comega_{10}}{5}$ are also
diagonal, $x^2+2y^2+3z^2+6w^2$ and $x^2+2y^2+5z^2+10w^2$, and they
are universal.

Two Hermitian lattices $\qf{1,1}$ and $\qf{1,2}$ over $\Q{-3}$ are
associated with quadratic forms $x^2+xy+y^2+z^2+zw+w^2$ and
$x^2+xy+y^2+2z^2+2zw+2w^2$. They represent universal quadratic forms
$x^2+z^2+3y^2+3w^2$ and $x^2+2z^2+3y^2+6w^2$, respectively. Also,
$\qf{1}\perp\binlattice{2}{-1+\omega_5}{-1+\comega_5}{3}$ over
$\Q{-5}$ contains a universal lattice $\qf{1,2}$.

The quadratic forms associated to $\qf{1,1}$ over $\Q{-7}$ and
$\qf{1,2}$ over $\Q{-11}$ lie in one class genera as listed in
\cite{Nipp}.

\begin{align*}
\Q{-7}:\,
    &\qf{1,1}
        \text{ corresponds to } x^2+y^2+2z^2+2w^2+xz+yw \\
\Q{-11}:\,
    &\qf{1,2}
        \text{ corresponds to } x^2+2y^2+3z^2+6w^2+xz+2yw \\
\end{align*}

Thus nine lattices remain:
\begin{align*}
\Q{-7}:\,
    &\qf{1,2}
        \text{ corresponds to } f_{7,2} = x^2+2y^2+2z^2+4w^2 +xy+2zw, \\
    &\qf{1,3}
        \text{ corresponds to } f_{7,3} = x^2+2y^2+3z^2+6w^2 +xy+3zw, \\
\Q{-11}:\,
    &\qf{1,1}
        \text{ corresponds to } f_{11} = x^2+y^2+3z^2+3w^2 +xz+yw, \\
\Q{-15}:\,
    &\qf{1}\perp\binlattice{2}{\omega_{15}}{\comega_{15}}{2} \\%
    &\phantom{\qf{1,1}}
        \text{ corresponds to } f_{15} = x^2+2y^2+2z^2+4w^2 +xw+yz, \\
\Q{-19}:\,
    &\qf{1,2}
        \text{ corresponds to } f_{19} = x^2+2y^2+5z^2+10w^2 +xz+2yw, \\
\Q{-23}:\,
    &\qf{1}\perp\binlattice{2}{\omega_{23}}{\comega_{23}}{3},
        \qf{1}\perp\binlattice{2}{-1+\omega_{23}}{-1+\comega_{23}}3 \\%
    &\phantom{\qf{1,1}}
        \text{ correspond to } f_{23} = x^2+2y^2+3z^2+6w^2 +xw+yz, \\
\Q{-31}:\,
    &\qf{1}\perp\binlattice{2}{\omega_{31}}{\comega_{31}}{4},
        \qf{1}\perp\binlattice{2}{-1+\omega_{31}}{-1+\comega_{31}}{4} \\%
    &\phantom{\qf{1,1}}
        \text{ correspond to } f_{31} = x^2+2y^2+4z^2+8w^2 +xw+yz. \\
\end{align*}
Note that the quadratic forms associated to each two lattices over
$\Q{-23}$ and $\Q{-31}$ coincide. So we will show the universalities
of the seven quaternary quadratic forms. The key idea is to find a
genus whose classes are all represented by the associated quaternary
quadratic form. To do this we use the Brandt-Intrau-Schiemann tables
\cite{Brandt-Intrau-Schiemann} via computer search.

\subsection{Universality of $f_{7,2}$}

Let $g = x^2+2y^2+2z^2 + xy$. Then $f_{7,2}$ represents $g(x,y,z)+
2\cdot7w^2$. Note that $g$ represents two forms $g_1 =
x^2+9y^2+15z^2+6yz$ and $g_2 = 3x^2+6y^2+7z^2$ which constitute a
genus \cite{Brandt-Intrau-Schiemann}. So $g$ represents all positive
integers which is represented by the genus of $g_1$ and $g_2$. That
is, $g$ represents $n$ unless $n \equiv 2 \pmod{3}$ or $7 | n$.
Since $7 = (-1+2\omega)(-1+2\comega)$, $7^s\cdot3t$ with $7 \nmid t$
and $7^s(3t+1)$ with $7 \nmid (3t+1)$ are represented by $f_{7,2}$.
If $n = 3t+2$ with $7 \nmid (3t+2)$, then $n-2\cdot7 \nequiv 2
\pmod{3}$ and $7 \nmid (n-2\cdot7)$. Thus $f_{7,2}$ represents
$7^s(3t+2)$. Since $f_{7,2}$ represents all positive integers up to
$2\cdot7$, $f_{7,2}$ is universal.

\subsection{Universality of $f_{7,3}$}

The quadratic form $f_{7,3}$ represents two ternary forms $g_1 =
x^2+3y^2+7z^2$ and $g_2 = 2x^2+3y^2+4z^2 + 2xz$ which constitute a
genus \cite{Brandt-Intrau-Schiemann}. So $f_{7,3}$ represents all
positive integers not divisible by $3$. Since $\qf{1,3}$ represents
$3\qf{1,3}=\qf{3,3^2}$, $f_{7,3}$ represents all positive integers.

\subsection{Universality of $f_{11}$}

The quadratic form $f_{11}$ represents two ternary forms $g_1 =
x^2+y^2+11z^2$ and $g_2 = x^2+3y^2+4z^2 + 2yz$ which constitute a
genus \cite{Brandt-Intrau-Schiemann}. Thus $f_{11}$ represents all
positive integers not divisible by $11$. Since $11 =
(-1+2\omega)(-1+2\comega)$, $f_{11}$ represents all positive
integers.

\subsection{Universality of $f_{15}$}

The quadratic form $f_{15}$ represents $g_1 = x^2+2y^2+8z^2 +2yz$
and $g_2 = x^2+3y^2+5z^2$ which constitute a genus
\cite{Brandt-Intrau-Schiemann}. So $f_{15}$ represents all positive
integers not divisible by $5$. Since $f_{15}(2y+3z,x+3w,x-2w,y-z) =
5f_{15}(x,y,z,w)$, $f_{15}$ represents all positive integers.

\subsection{Universality of $f_{19}$}

Let $g = x^2+2y^2+5z^2 +xz$. Then $f_{19}$ represents
$g(x,y,z)+2\cdot19w^2$. Note that $g$ represents two forms $g_1 =
2x^2+5y^2+25z^2 +5yz$ and $g_2 = 3x^2+7y^2+13z^2 +3yz+xz+3xy$ which
constitute a genus \cite{Brandt-Intrau-Schiemann}. So $g$ represents
all positive integers $n$ unless $n \equiv 1, 4 \pmod{5}$ or $2 |
n$. Since $\qf{1,2}$ represents $2\qf{1,2}$, $f_{19}$ represents all
positive integers $n = 2^s(5t+k)$ with $k=2,3$ and $2 \nmid (5t+k)$.
Suppose $n = 5t+1$ with $t \geq 8$ is odd. Then $n - 2\cdot19 =
5(t-8)+3$ is represented by $g$. If $n = 5t+4$ with $t \geq 30$ is
odd, then $n - 2\cdot19\cdot2^2 = 5(t-30)+2$ is represented by $g$.
Since $\qf{1,2}$ represents all positive integers up to
$5\cdot30+4$, $f_{19}$ is universal.

\subsection{Universality of $f_{23}$}

Let $g = x^2+2y^2+3z^2 +yz$. Then $f_{23}$ represents
$g(x,y,z)+23w^2$. Note that $g$ represents two forms $g_1 =
x^2+8y^2+12z^2 +4yz$ and $g_2 = 4x^2+4y^2+9z^2+ 4yz+4xz+4xy$ which
constitute a genus \cite{Brandt-Intrau-Schiemann}. Thus $g$
represents all positive integers $n$ unless $n \equiv 2,3 \pmod{4}$
or $23 | n$. Since $f_{23}(y,2z,2x,w) = 2f_{23}(x,y,z,w)$ and $23 =
(-1+2\omega)(-1+2\comega)$, $f_{23}$ represents $n = 2^r \cdot 23^s
(4t+1)$ with $23 \nmid (4t+1)$. Suppose that $n=4t+3$ is not
divisible by $23$. Then $n-23 = 4(t-5)$ is represented by $g$ since
$23 \nmid (t-5)$. Since $f_{23}$ represents all positive integers up
to $4\cdot5+3$, $f_{23}$ is universal.

\subsection{Universality of $f_{31}$}

Let $g = x^2+2y^2+4z^2 +yz$. Then $f_{31}$ represents
$g(x,y,z)+31w^2$. Note that $g$ represents three forms $g_1 =
x^2+4y^2+32z^2 +4yz$, $g_2 = x^2+8y^2+16z^2+ 4yz$ and $g_3 =
4x^2+5y^2+8z^2 +4yz+4xz$ which constitute a genus
\cite{Brandt-Intrau-Schiemann}. This genus represents all positive
integers $n$ unless $n \equiv 2,3 \pmod{4}$ or $31 | n$. Since
$f_{31}(2y,w,x,2z) = 2f_{31}(x,y,z,w)$ and $31 =
(-1+2\omega)(-1+2\comega)$, $f_{31}$ represents $n = 2^r \cdot 31^s
(4t+1)$ with $31 \nmid (4t+1)$. Suppose that $n=4t+3$ is not
divisible by $31$. Then $n-31 = 4(t-7)$ is represented by $g$ since
$31 \nmid (t-7)$. Since $f_{31}$ represents all positive integers up
to $4\cdot7+3$, $f_{31}$ is universal.

\end{document}